\documentclass[11pt,reqno]{article}
\usepackage{amsmath,mathrsfs,graphicx,amssymb,amsfonts,color}

\font\bg=cmbx10 scaled\magstep1

\font\small=cmr8

\newtheorem{newlemma}{{\bf Lemma}}

\newenvironment{lema}{\begin{newlemma}{\hspace{-0.5
em}{\bf.}}}{\end{newlemma}}

\newtheorem{newteorem}{{\bf Theorem}}

\newenvironment{teorem}{\begin{newteorem}{\hspace{-0.5
em}{\bf.}}}{\end{newteorem}}

\newtheorem{newkorolari}{{\bf Corollary}}

\newenvironment{korolari}{\begin{newkorolari}{\hspace{-0.5
em}{\bf.}}}{\end{newkorolari}}

\newtheorem{newdefine}{{\bf Definition}}

\newenvironment{define}{\begin{newdefine}{\hspace{-0.5
em}{\bf.}}}{\end{newdefine}}

\newtheorem{newquestion}{{\bf Question}}

\newtheorem{newkonjek}{{\bf Conjecture}}

\newtheorem{newexample}{{\bf Example}}

\begin{document}
\tolerance=10000
\baselineskip18truept
\newbox\thebox
\def\boxit#1{\vbox{\hrule\hbox{\vrule\kern0pt
\vbox{\kern0pt#1\kern0pt}\kern0pt\vrule}\hrule}}
\def\qed{\lower0.1cm\hbox{\noindent \boxit{\copy\thebox}}\bigskip}
\def\ss{\smallskip}
\def\ms{\medskip}
\def\bs{\bigskip}
\def\c{\centerline}
\def\nt{\noindent}
\def\ul{\underline}
\def\ol{\overline}
\def\lc{\lceil}
\def\rc{\rceil}
\def\lf{\lfloor}
\def\rf{\rfloor}
\def\ov{\over}
\def\t{\tau}
\def\th{\theta}
\def\k{\kappa}
\def\l{\lambda}
\def\L{\Lambda}
\def\g{\gamma}
\def\d{\delta}
\def\D{\Delta}
\def\e{\epsilon}
\def\lg{\langle}
\def\rg{\rangle}
\def\p{\prime}
\def\sg{\sigma}
\def\ch{\choose}

\newcommand{\ben}{\begin{enumerate}}
\newcommand{\een}{\end{enumerate}}
\newcommand{\bit}{\begin{itemize}}
\newcommand{\eit}{\end{itemize}}
\newcommand{\bea}{\begin{eqnarray*}}
\newcommand{\eea}{\end{eqnarray*}}
\newcommand{\bear}{\begin{eqnarray}}
\newcommand{\eear}{\end{eqnarray}}

\centerline{\Large \bf Domination polynomials of $k$-tree related graphs   }
\vspace{.3cm}

\bigskip

\bs

\baselineskip12truept
\centerline{ S. Jahari  and S. Alikhani$^{}${}\footnote{\baselineskip12truept\it\small
Corresponding author. E-mail: alikhani@yazd.ac.ir} }
\baselineskip20truept
\centerline{\it Department of Mathematics, Yazd University}
\vskip-8truept
\centerline{\it  89195-741, Yazd, Iran}

\vskip-0.2truecm
\nt\rule{16cm}{0.1mm}

\nt{\bg ABSTRACT}
\medskip

\baselineskip14truept

\nt{ Let $G$ be a simple graph of order $n$.
The domination polynomial of $G$ is the polynomial
$D(G, x)=\sum_{i=\gamma(G)}^{n} d(G,i) x^{i}$,
where $d(G,i)$ is the number of dominating sets of $G$ of size $i$ and
$\gamma(G)$ is the domination number of $G$.
In this paper we study  the domination polynomials of several classes of $k$-tree related graphs. Also,  we present  families
of these kind of graphs, whose domination polynomial have no nonzero real roots.}
\ms

\nt{\bf Mathematics Subject Classification:} {\small 05C60.}

\nt\rule{16cm}{0.1mm}

\baselineskip20truept

\section{Introduction}

\nt Throughout this paper we will consider only simple graphs. Let $G=(V,E)$ be a simple graph.
 For $F\subseteq V(G)$ we use $<F>$ for the subgraph induced by $F$. For any vertex $v\in V(G)$, the {\it open neighborhood} of $v$ is the
set $N(v)=\{u \in V (G) | \{u, v\}\in E(G)\}$ and the {\it closed neighborhood} of $v$
is the set $N[v]=N(v)\cup \{v\}$. For a set $S\subseteq V(G)$, the open
neighborhood of $S$ is $N(S)=\bigcup_{v\in S} N(v)$ and the closed neighborhood of $S$
is $N[S]=N(S)\cup S$. For every vertex $v\in V (G)$, the {\it degree} of $v$ is the number of
edges incident with $v$ and is denoted by $d_G(v)=|N(v)|$. Let $d_i, ~1\leq i \leq n$, be the degrees of the vertices $v_i$ of a graph in any order. The sequence $\{d_i\}^n_1$ is called the {\it degree sequence} of the graph.  A {\it clique} in a graph is a subset of its vertices such that every two vertices in the subset are connected by an edge. We use $K_n, P_n, C_n$ and $S_{1,n-1}$ for a clique, a path, a cycle
and a star, all of order $n$, respectively.

\nt A set $S\subseteq V(G)$ is a {\it dominating set} if $N[S]=V$ or equivalently,
every vertex in $V(G)\backslash S$ is adjacent to at least one vertex in $S$.
The {\it domination number} $\gamma(G)$ is the minimum cardinality of a dominating set in $G$. A dominating set
with cardinality $\gamma(G)$ is  called a {\it $\gamma$-set}.
For a detailed treatment of these parameters, the reader is referred to~\cite{domination}.
Let ${\cal D}(G,i)$ be the family of dominating sets of a graph $G$ with cardinality $i$ and
let $d(G,i)=|{\cal D}(G,i)|$.
The {\it domination polynomial} $D(G,x)$ of $G$ is defined as
$D(G,x)=\sum_{ i=\gamma(G)}^{|V(G)|} d(G,i) x^{i}$,
where $\gamma(G)$ is the domination number of $G$ (see \cite{euro,saeid1}).
 Thus $D(G, x)$ is the generating polynomial for the number of dominating sets of $G$
of each cardinality.
A root of $D(G, x)$ is called a {\it domination root} of $G$.

\ms

\nt In \cite{Kot} it is
shown that computing the domination polynomial $D(G, x)$ of a graph $G$ is NP-hard and
some examples for graphs for which $D(G, x)$ can be computed efficiently are given.
 The vertex contraction $G/u$ of a graph $G$ by a vertex $u$ is the operation under
which all vertices in $N(u)$ are joined to each other and then $u$ is deleted (see \cite{Wal}).
  The following theorem is  useful for finding the recurrence relations for the  domination polynomials  of  arbitrary graphs.

\begin{teorem}\label{theorem1}{\rm \cite{saeid2,Kot}}
Let $G$ be a graph. For any vertex $u$ in $G$ we have
\[
D(G, x) = xD(G/u, x) + D(G - u, x) + xD(G - N[u], x) - (1 + x)p_u(G, x),
\]
where $p_u(G, x)$ is the polynomial counting the dominating sets of $G - u$ which do not contain any
vertex of $N(u)$ in $G$.
\end{teorem}

\nt Using Theorem \ref{theorem1} we are able to obtain an easier formula for a graph with at lease one vertex of degree 1. Since every
tree has at least two vertices of degree 1, so we can use the following recurrence to obtain the domination polynomials of trees.

\begin{korolari}{\rm \cite{JME,Kot}}
 Let $G = (V,E)$ be a graph, $v$ be a vertex of degree $1$ in $G$ and let $u$ be its
neighbor. Then
\[
D(G, x) = xD(G/u, x) + D(G - u - v, x) + D(G - N[u], x).
\]
\end{korolari}

\nt If $G_1$ and $G_2$ are disjoint graphs
of orders $n_1$ and $n_2$  respectively, then $D(G_1 \cup G_2, x) = D(G_1, x)D(G_2, x)$ and
\[
D(G_1 + G_2,x)=\Big((1+x)^{n_1}-1\Big)\Big((1+x)^{n_2}-1\Big)+D(G_1,x)+D(G_2,x),
\]
where $G_1+ G_2$ is the join of $G_1$ and $G_2$, formed from $G_1 \cup G_2$ by adding in all edges
between a vertex of $G_1$ and a vertex of $G_2$ (see \rm\cite{euro}).

\nt The domination polynomials of trees, aside from path graph have not been studied and there is no study for coefficients
of $D(T,x)$ for trees $T$ with $n$ vertices. $k$-trees are generalization of tree which consider in this paper.
Actually similar to  \cite{song}, in this paper  we consider  $k$-tree related graphs and study their domination polynomials.
Study of the roots of domination polynomials  is a interesting (\cite{SSM,brown}). One of the problem in domination roots is classification and finding graphs with no nonzero real roots. In this paper we present some families related to $k$-trees which have this property.


\ms
\nt In Section 2, we study the  domination polynomials for some $k$-tree related graphs.
In Section 3, we present some families
of these kind of graphs whose domination polynomials have no nonzero real roots.

\section{Domination polynomials of $k$-tree related graphs }

In this section we study the  domination polynomials for some $k$-tree related graphs. The class of $k$-trees is a very important
subclass of triangulated graphs. Harary and Palmer \rm\cite{harary} first introduced $2$-trees in 1968. Beineke and Pippert \rm\cite{bein} gave the definition of a $k$-tree in 1969. In the literature on $k$-trees, there are interesting applications to the study of computational complexity.

\begin{define} For a positive integer $k$, a $k$-tree, denoted by $T^k_n$, is defined recursively as follows: The smallest $k$-tree is the $k$-clique $K_k$. If $G$ is a $k$-tree with $n\geq k$ vertices and a new vertex $v$ of degree $k$ is added and joined to the vertices of a $k$-clique in $G$, then the larger graph is a $k$-tree with $n +1$ vertices.
\end{define}

\begin{figure}[!h]
\hspace{2.9cm}
\includegraphics[width=9.6cm,height=2.4cm]{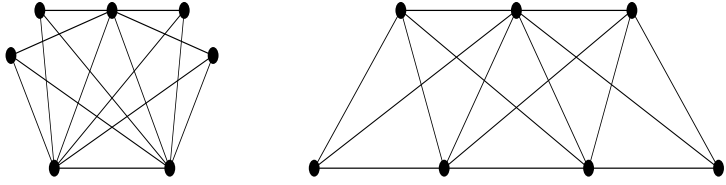}
\caption{\label{Figure1} The 3-star and  3-path on 7 vertices, respectively. }
\end{figure}

\nt An {\it independent set} in a graph $G$ is a set of pairwise non-adjacent vertices.
\begin{define}
Let $K_k$ be a $k$-clique and $S$ be an independent set of $n- k$ vertices. A $(k, n)$-star, denoted by $S_{k,n-k}$, is defined
as $S_{k,n-k}=K_k +S$.
\end{define}
\begin{define}
 A $(k, n)$-path, denoted by $P^k_n$, begins with $k$-clique on $\{v_1, v_2, \dots, v_k\}$. For $i= k+ 1$ to $n$, let vertex $v_i$ be adjacent to vertices $\{v_{i-1}, v_{i-2}, \dots, v_{i-k}\}$ only. (see Figure \ref{Figure1}).
\end{define}

\nt A helpful characteristic of the $(k, n)$-path  $P^k_n$ is that we may order the vertices $v_1, v_2, \dots, v_n$ such that  $P^k_n \setminus \{v_1, \dots, v_i\}$ is a $k$-path on $n - i$ vertices for $1 \leq i \leq n - k - 1$, such a vertex ordering is referred to as a presentation.
\begin{define}
 A $(k, n)$-cycle, denoted by $C^k_n$, consists of a $(k, n)$-path on $\{v_1, v_2, \dots, v_n\}$ defined as above and an edge joining $v_1$ to $v_n$, where $n\geq k + 2$.
\end{define}
\begin{define}
If $G$ is a $(k, n)$-cycle of order $n$ and $v$ is a vertex not in $G$, then $G+v$ is called a $(k, n)$-wheel and denoted by $W^k_n$.
\end{define}

\nt Notice that $P^1_n,~ C^1_n,~ W^1_n$ and $S_{1,n-1}$ are just the standard path, cycle, wheel and star, respectively.
It follows easily from the domination polynomial of join of two graphs that,  for the
star graph $K_{1,n-1}$, we have $D(K_{1,n},x) = x(1 + x)^{n-1} + x^{n-1}.$
\nt The following are  recurrences  for the domination polynomials of  paths and cycles (\cite{peng}).
\begin{teorem} For the natural number $n\geq 3$, \\
$(i)~D(P_{n+1}, x) = x\big(D(P_n, x) + D(P_{n-1}, x) + D(P_{n-2}, x)\big),$ \\
where $D(P_1, x) = x,~ D(P_2, x) = x^2 + 2x$ and $D(P_3, x) = x^3 + 3x^2 +x $. \\
$(ii)~D(C_{n+1}, x) = x\big(D(C_n, x) + D(C_{n-1}, x) + D(C_{n-2}, x)\big),$ \\
where $D(C_1, x) = x,~ D(C_2, x) = x^2 + 2x$ and $D(C_3, x) = x^3 + 3x^2 +3x $.
\end{teorem}

\nt Note that both $(k, n)$-cycles and $(k, n)$-wheels are not $k$-trees. But they are closely related to $k$-trees. We begin by a
simple lemma was proven in \rm\cite{song} as Proposition 2:

\begin{lema}
 For any $k$-tree $T^k_n$, $| E(T^k_n) | = kn -\frac{1}{2}k(k + 1) = \frac{2nk-k^2-k}{2}$.
\end{lema}

\nt The {\it independence number} is the size of a maximum independent set in the graph and denoted by $\alpha(G)$. The following Lemma gives independence numbers  for $k$-tree related graphs:

\begin{lema}\label{i}\rm\cite{song} For each natural number $k\leq n$, we have \\
$(i)~ \alpha(P^k_n) = \lfloor \frac{n+k}{k +1} \rfloor$.\\
$(ii)~ \alpha(C^k_n) = \lfloor \frac{n+k-1}{k +1} \rfloor$.\\
$(iii)~ \alpha(W^k_n) = \alpha(C^k_n)$.\\
$(iv)~ \alpha(S_{k,n-k}) = n - k$.
\end{lema}

\nt Now, we present the following domination numbers  for $k$-tree related graphs:

\begin{teorem} For each natural number $k\leq n$, we have \\
$(i)~ \gamma(P^k_n) = \lceil \frac{n}{2k +1} \rceil$.\\
$(ii)~ \gamma(C^k_n) = \lceil \frac{n}{2k +1} \rceil$.\\
$(iii)~ \gamma(W^k_n) = 1$.\\
$(iv)~ \gamma(S_{k,n-k}) = 1$.
\end{teorem}
\nt{\bf Proof.}
$(i)$ Since $k\leq n$,  by the definition of $P^k_n$, the degree sequence in this graph is
$$\{k, k+1, \dots, 2k-1, 2k, \dots, 2k, 2k-1, \dots, k-1, k\},$$
 we have   $\gamma(P^k_n) = 1$ for  $n\leq 2k+1$ and   $\gamma(P^k_n) =2$ for  $~2k+2\leq n\leq 4k+2$. Thus $(i)$ holds for $k \leq n \leq 4k+2$. Now assume $n \geq 4k+3$. We use induction on $n$. Since any $\gamma$-set of $P^k_n$ contains only one vertex of the $\{v_n, v_{n-1}, \dots, v_{n-k}, \dots, v_{n-2k-1}\}$ and $P^k_n - \{v_n, v_{n-1}, \dots, v_{n-k}, \dots, v_{n-2k-1}\}$ is a $k$-path with $n - 2k - 1$ vertices, by induction, $\gamma(P^k_n) = 1 + \gamma(P^k_{n-2k-1}) = 1 + \lceil \frac{n-2k-1}{2k +1} \rceil = \lceil \frac{n}{2k +1} \rceil$. Hence $(i)$ holds.

\nt $(ii)$ Since $k\leq n$,  by the definition of $C^k_n$, the degree sequence in this graph is
$$\{k+1, k+1, k+2, \dots, 2k-1, 2k, \dots, 2k, 2k-1, \dots, k-1, k+1\},$$
 we have   $\gamma(C^k_n) = 1$ for  $n\leq 2k+1$ and   $\gamma(C^k_n) =2$ for   $~2k+2\leq n\leq 4k+2$. Thus $(ii)$ holds for $k \leq n \leq 4k+2$. Now assume $n \geq 4k+3$ and use induction on $n$. Since any $\gamma$-set of $C^k_n$ contains only one vertex of the $\{v_n, v_{n-1}, \dots, v_{n-k}, \dots, v_{n-2k-1}\}$ and $C^k_n - \{v_n, v_{n-1}, \dots, v_{n-k}, \dots, v_{n-2k-1}\}$ is a $k$-path with $n - 2k - 1$ vertices, by induction, $\gamma(C^k_n) = 1 + \gamma(P^k_{n-2k-1}) = 1 + \lceil \frac{n-2k-1}{2k +1} \rceil = \lceil \frac{n}{2k +1} \rceil$. Hence $(ii)$ holds.

\nt $(iii)$ Since the $(k, n)$-wheel $W^k_n$, has a vertex $v$ of degree  $n-1$, so   $(iii)$ holds. 

\nt $(iv)$ Since the $k$-star graph $S_{k,n-k}$, has $k$ vertices of degree  $n-1$, so   $(iv)$ holds.$\quad\Box$

\nt  The  following theorem gives a recurrence formula for the domination polynomial of $(k, n)$-path graphs.

\begin{teorem}\label{p} If $ n \leq k+1$, then $D(P^k_n, x) =D(K_n,x)$.
For every $k+2\leq n$,
\[D(P^k_n, x) = (1+x)D(P^k_{n-1}, x) +x D(P^k_{n-k-1}, x) - (1+x)p_u(P^k_{n}, x),\]
where
$
p_u(P^k_{n}, x)=\left\{
\begin{array}{lr}
{ x(1+x)^{n-k-2}};&
\quad\mbox{$k+2\leq n\leq 2k+2$,}\\[15pt]
{x((1+x)^{n-k-2}-(1+x)^{n-2k-3})};&
\quad\mbox{$2k+3\leq n\leq 2k+6$,}\\[15pt]
{p_u(P^k_{n}, x)};&\quad\mbox{$2k+7\leq n$.}
\end{array}
\right.
$
\end{teorem}
\nt{\bf Proof.} If $n \leq k+1$, then $P^k_n\cong K_n$.
For every $k+2\leq n$,   we use Theorem \ref{theorem1} for the last vertex of $P^k_n$ and since (by the definition of $P^k_n$) the first $k+1$ and the last $k+1$ vertices form two cliques, we have $P^k_{n}/u \cong P^k_{n} - u$. It is clear that
$D(P^k_n - N[u], x) = D(P^k_{n-k-1}, x)$. Obviously $p_u(P^k_{n}, x)$ is the polynomial counting
the dominating sets of $P^k_{n-k-1}$ contains the vertex  $v_{n-k-1}$, but finding  this polynomial
 involve complex calculations. 
We brought this polynomial for $n\leq 2k+6$ in this theorem. Therefore we have the result.$\quad\Box$

\nt In general, finding the domination polynomial of a graph is a very difficult problem.
In \rm\cite{Kot} Kotek et al. showed that there exist recurrence relations for the domination polynomial
which allow for efficient schemes to compute the polynomial for some types of graphs.
Consider $(k, n)$-cycle graphs, If $n \leq k+2$, then $C^k_n\cong K_n$. Consequently in this case $D(C^k_n, x) =D(K_n,x)$.
For every $k+3\leq n$, until now all attempts to find formulas
for $D(C^k_n, x)$ failed.  

  \nt  The  following theorem gives a formula for the domination polynomial of $(k, n)$-wheel graphs.

\begin{teorem}\label{w} For a $(k,n)$-wheel $W^k_n$, we have
\[D(W^k_n, x) =x(1+x)^{n-1} + D(C^k_n, x).\]
\end{teorem}
\nt{\bf Proof.} Since $W^k_n = C^k_n+ K_1$, then
\begin{eqnarray*}
D(W^k_n, x)&=&((1+x)-1)((1+x)^{n-1}-1) + x + D(C^k_n, x)\\
&=&x(1+x)^{n-1} + D(C^k_n, x).\quad\quad\Box
\end{eqnarray*}

\nt  The  following theorem gives a formula for the domination polynomial of $k$-star graphs, which is concluded of the fact,
$k$-star graph is the join of complete graph $K_k$ and independent set $S$ (empty graph $O_{n- k}$).

\begin{teorem}\label{theorem}
For every $k\in \mathbb{N}$ and $n> k$,
\[D(S_{k,n-k}, x) = (1+x)^{n-k}((1+x)^k -1) +x^{n-k}. \]
\end{teorem}
\nt{\bf Proof.} Let $S_{k,n-k}$ be the  $k$-star graph with vertex set $V(S_{k,n-k}) = \{v_1, v_2, \cdots , v_k, v_{k+1},  \cdots , v_n\}$. It suffices to show that every dominating set of
size $j = 1,\ldots,n$ is accounted for exactly once in the above statement.
Clearly every non-empty subset of $\{v_1, v_2, \cdots , v_k\}$ is dominating set of $k$-star graphs. These sets can be extended with any number of vertices $\{v_{k+1},  \cdots , v_n\}$. Also obviously  the set $\{v_{k+1},  \cdots , v_n\}$ is a dominating set of $k$-star graphs. It is easy to see that there is no another method to make a dominating set for $k$-star graphs. Therefore we have the result.$\quad\Box$

\nt The value of a graph polynomial at a specific point can give
sometimes a surprising information about the structure of the graph \cite{gcom,levit}. The following simple results give
the domination polynomial  of $k$-tree related graphs at  $-1$.

\begin{korolari} For each natural number $k\leq n$, the following hold: \\
$(i)~ D(P^k_n, -1) = (-1)^{\alpha(P^k_n) }$.\\
$(ii)~ D(W^k_n, -1) = D(C^k_n, -1)$.\\
$(iii)~ D(S_{k,n-k, -1}) =(-1)^{\alpha(S_{k,n-k}) }$.
\end{korolari}
\nt{\bf Proof.} $(i)~$ Using domination polynomial of $(k, n)$-path in Theorem \ref{p}, for $k\leq n \leq k+1$, $D(P^k_n, -1) = D(K_n,-1) = -1$. For every $k+2\leq n$, $D(P^k_n, -1) = - D(P^k_{n-k-1}, -1)$. Obviously, in the first case  $\lfloor \frac{n+k}{k +1} \rfloor = 1$. Thus $(i)$ holds for $ n \leq k+1$. Now assume $k+2\leq n$. We use induction on $n$. Suppose that the statement is true for every $k$-path with $n - k$ vertices, by induction and Lemma \ref{i},
 \begin{eqnarray*}
 D(P^k_n, -1)&=& - D(P^k_{n-k-1}, -1)\\
 &=& - (-1)^{\alpha(P^k_{n-k-1}) } = - (-1)^{\lfloor \frac{n-1}{k +1} \rfloor}\\
 &=& (-1)^{1+\lfloor \frac{n-1}{k +1} \rfloor} = (-1)^{\lfloor \frac{n+k}{k +1} \rfloor}.
 \end{eqnarray*}
 Hence $(i)$ holds.\\
$(ii)~$ Follows from  Theorem \ref{w}.\\
$(iii)~$ Follows from   Theorem \ref{theorem} and Lemma \ref{i}.$\quad\Box$

\section{Some families of graphs with no nonzero real domination roots}

\nt In \rm\cite{SSM} authors asked that which graphs have no nonzero real domination roots?

\nt In this section we would like to obtain more results related to this problem.  We need some preliminaries.

\nt For two graphs $G = (V,E)$ and $H=(W,F)$, the corona $G\circ H$ is the graph arising from the
disjoint union of $G$ with $| V |$ copies of $H$, by adding edges between
the $i$th vertex of $G$ and all vertices of $i$th copy of $H$ \cite{Fruc}. It is easy to see that the corona operation of two graphs does not have the commutative property.

\nt We need the following theorem which is for computation of domination
polynomial of corona products of two graphs.

\begin{teorem}\label{theorem7}{\rm \cite{Oper,Kot}}
Let $G = (V,E)$ and $H=(W,F)$ be nonempty graphs of order $n$ and $m$, respectively. Then
\begin{eqnarray*}
D(G\circ H,x) = (x(1 + x)^m + D(H, x))^n.
\end{eqnarray*}
\end{teorem}

\nt A $k$-star, $S_{k,n-k}$, has vertex set $\{v_1, \cdots , v_n\}$ where $<\{v_1, v_2, \cdots , v_k\}>  \cong K_k$ and $N(v_i) = \{v_1, \cdots, v_k\}$ for $k + 1 \leq i \leq n$.


\nt Here we will discuss roots of domination polynomial of $k$-star graphs.

\begin{teorem}
\begin{enumerate}
\item[(i)]
For odd  natural  $n$ and even  natural  $k$, no nonzero real numbers is  domination root of $S_{k,n-k}$.

\item[(ii)]
For even  natural  $n$ and even  natural  $k$, there is exactly one nonzero real domination root of $S_{k,n-k}$.

\end{enumerate}
\end{teorem}
\nt{\bf Proof.} By Theorem \ref{theorem}, for every $n > k$, $D(S_{k,n-k}, x) = (1+x)^{n-k}((1+x)^k -1) +x^{n-k}$. If
$D(S_{k,n-k},x)=0$, then we have
\[
(1+x)^{n-k}((1+x)^k -1) = - x^{n-k}.
\]
Now we are ready to prove two cases of this theorem:
\begin{enumerate}
 \item[(i)]
First suppose that $x\geq 0$. Obviously the above equality is true
just for real number 0, since for nonzero real number the left
side of equality is positive but the right side is negative. Now
suppose that $x< -1$. In this case the left side is negative
and the right side $-x^{n-k}$ is greater
than $+1$, a contradiction. Finally we shall consider $-1<x<0$.
This case is similar to the second case when we substitute $x$
with $\frac{1}{x}$.
\item[(ii)]
First suppose that $x\geq 0$ and $x< -1$. Obviously the above equality is true
just for real number 0, since for nonzero real number the left
side of equality is positive but the right side is negative. Now
suppose that $-1<x<0$. This equation has
only one real root in $(-1, 0)$.$\quad\Box$
\end{enumerate}

\begin{figure}[h]
\hspace{4.cm}
\includegraphics[width=6.3cm,height=5.cm]{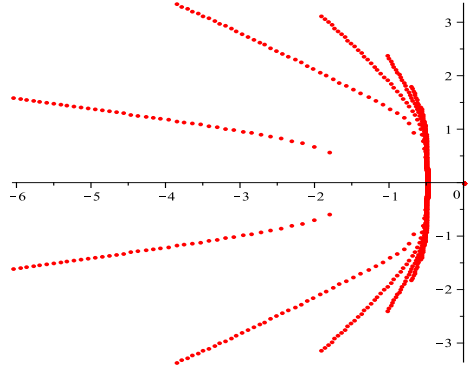}
\caption{\label{figure2} Domination roots of graphs $S_{4,n-4}$ for $5 \leq n \leq 44$.}
\end{figure}

\nt{\bf Remark.} Using Maple we have shown the domination roots
 of $S_{4,n-4}$ for $5\leq n\leq 44$ in  Figure \ref{figure2}.

\nt Here we construct a sequence of graphs, which their  domination roots are  the same as the  domination roots of the $k+1$-star graphs.

\begin{teorem} The domination roots of every graph $H$ in the family
$$\{G\circ S_{k, n-k}, (G\circ S_{k, n-k})\circ S_{k, n-k}, ((G\circ S_{k, n-k})\circ S_{k, n-k})\circ S_{k, n-k},\cdots \}$$ have the same behavior as the domination roots of $k+1$-star graphs.
\end{teorem}
\nt{\bf Proof.}
 By theorem \ref{theorem7} we can deduce that for each arbitrary graph $G$,
\begin{eqnarray*}
D(G\circ S_{k, n-k},x)&=&\Big(x(1+x)^n + (1+x)^{n-k}((1+x)^k -1) +x^{n-k}\Big)^{|V(G)|}\\
 &=&\Big((1+x)^{n-k}((1+x)^{k+1} -1) +x^{n-k}\Big)^{|V(G)|}\\
 &=& \Big(D(S_{k+1, n-k},x))\Big)^{|V(G)|}.
\end{eqnarray*}
 Therefore we have the result.$\quad\Box$

\end{document}